\input amssym.def
\input amssym
\magnification=1200
\parindent0pt
\hsize=16 true cm
\baselineskip=13  pt plus .2pt
$ $

\def\C{\Bbb C}
\def\R{\Bbb R}
\def\Z{\Bbb Z}

\def\A{\Bbb A}
\def\S{\Bbb S}

\centerline {\bf On topological actions of finite, non-standard groups on
spheres}

\bigskip

\centerline {Bruno P. Zimmermann}

\bigskip

\centerline {Universit\`a degli Studi di Trieste}

\centerline {Dipartimento di Matematica e Geoscienze}

\centerline {34127 Trieste, Italy}

\bigskip \bigskip

Abstract.  {\sl The standard actions of finite groups on spheres $S^d$ are
linear actions, i.e. by finite subgroups of the orthogonal groups  ${\rm
O}(d+1)$. We prove that, in each dimension  $d>5$, there is a finite group $G$
which admits a faithful, topological action on a sphere $S^d$ but is not
isomorphic to a subgroup of ${\rm O}(d+1)$. The situation remains open for
smooth actions.}

\bigskip \bigskip

{\bf 1. Introduction}

\medskip

We are interested in faithful actions by homeomorphisms of finite groups on
spheres $S^d$. The standard actions on spheres are the linear or orthogonal
actions by the finite subgroups of the orthogonal group ${\rm O}(d+1)$. There
is a rich literature on smooth, non-linear actions of finite groups on spheres
(see the surveys [D], [Z1]), but  not much is known on the class of finite
groups which can occur for smooth or topological actions; in particular, not a
single example of a finite group seems to be known which admits a faithful,
smooth action on a sphere $S^d$ but does not admit a faithful, linear action on
$S^d$ (i.e., is not isomorphic to a subgroup of ${\rm O}(d+1)$). Concerning
topological actions,  our main result is the following.

\bigskip

{\bf Theorem.}  {\sl  For each dimension $d>5$, there is a finite group $G$
which admits a faithful, topological action on the sphere $S^d$ but is not
isomorphic to a subgroup of ${\rm O}(d+1)$.}

\bigskip

In fact, relying on an observation in [CKS], our methods would imply that in
each dimension $d \ge 5$ there are infinitely many such groups, see the remark
at the end of section 2. We note that the actions we construct are equivalent 
to simplicial actions but they are not locally linear and hence not equivalent
to smooth actions.

\medskip

The Theorem remains open in dimension 3.  In dimension 3, as a consequence of
the geometrization of finite group actions on 3-manifolds after Thurston and
Perelman (cf. [BLP],[DL]), every finite group which admits a faithful, smooth
action on $S^3$ is isomorphic (and even conjugate) to a subgroup of ${\rm
O}(4)$ (but this does not remain true for smooth actions on homology 3-spheres,
see [Z2] and [Z3, section 5] for a discussion). The major problem for
topological actions on the 3-sphere is then  the possible presence of wildly
embedded fixed point sets; for triangulable or PL actions, this phenomenon does
not occur, the actions are then locally linear, conjugate to smooth actions
and, by the geometrization, also to linear ones. On the other hand, almost
nothing seems to be known on the possible finite groups which can act on $S^3$
with wild fixed point sets.

\medskip

Concerning dimension 4, it is proved in [CKS]  (completing results in [MeZ1,2])
that every finite group which admits a smooth or locally linear,
orientation-preserving, faithful action on
$S^4$ or on a homology 4-sphere is isomorphic to a subgroup of ${\rm SO}(5)$,
but that this does not remain true for topological actions on $S^4$ (the
orientation-reversing analogue is still open in the smooth case, and again not
true in the topological case).

\medskip

Returning to arbitrary dimensions, we have the following:

\bigskip

{\bf Question.}   Is there a finite group $G$ which admits a faithful, smooth
action on a sphere $S^d$ but is not isomorphic to a subgroup of  ${\rm
O}(d+1)$? Does the Theorem remain true for faithful, smooth actions?

\bigskip

We note that, for certain classes of finite groups $G$, the minimal dimension
of a faithful, smooth action of $G$ on a homology sphere coincides with the
minimal dimension of a faithful, linear action on a sphere, e.g. for the linear
fractional groups ${\rm PSL}_2(p)$ ([GZ, Theorem 3]), for finite $p$-groups
([DoH]), for some classes of alternating groups and some other finite simple
groups (and in some cases also for purely topological actions, using Smith fixed
point theory and the Borel formula which hold in a purely topological
setting).

\medskip

The proof of the Theorem is based on the existence, due to Milgram [Mg], of a
finite group $Q$ (Milnor group) which admits a faithful, smooth action on a
homology 3-sphere $M^3$ but is not isomorphic to a subgroup of the orthogonal
group  ${\rm O}(4)$, and also on the double suspension theorem stating that the
double suspension or join  $M^3 * S^0 * S^0 \cong M^3 * S^1$  (see e.g. [Mu])
of a homology 3-sphere $M^3$ is homeomorphic to $S^5$ (see [Ca]); so this
allows to shift finite actions on homology 3-spheres to topological actions on
spheres in higher dimensions.

\bigskip

{\bf 2. Proof of the Theorem}

\medskip

By strong results of Milgram [Mg], there exists a finite group $Q$ which admits
a smooth, faithful action on a homology 3-sphere $M^3$ but is not isomorpic to a
subgroup of the orthogonal group ${\rm O}(4)$ (so $Q$ does not admit a
faithful, linear  action on $S^3$). The group $Q$ is a Milnor group
$Q(8a,b,c)$ ([Mn]), for relatively coprime odd positive integers $a, b$ and $c$
with  $a \ge 3$ and $b>c\ge 1$. Such a Milnor group $Q(8a,b,c)$ is an extension
of
$\Z_{abc}  \cong \Z_a \times \Z_b \times \Z_c$ by the quaternion group $Q(8) =
\{\pm 1,\pm i,\pm j,\pm k\}$ of order 8, where  $i, j$ and $k$ act trivially on
$\Z_a, \Z_b$ and $\Z_c$,  respectively, and in a dihedral way on the other two. 
By [Mn], $Q(8a,b,c)$ is not isomorphic to a subgroup of ${\rm O}(4)$ (see also
[Z3, section 3] for a discussion and  other references).

\medskip

By the double suspension theorem (see [Ca]), for $m \ge 1$ the join 
$$M^3 * S^m \;\; \cong \;\;  M^3 * S^1 * S^{m-2} \;\;  \cong \;\;  S^5 * S^{m-2}
\;\;  \cong \;\;  S^{m+4}$$ is homeomorphic to $S^{m+4}$. 

\medskip

The alternating group
$\A_n$ has a linear, faithful action on $\R^n$ by permutation of coordinates,
and also on $\R^{n-1}$  (in coordinates, on the subspace of $\R^n$ given by
$x_1+
\ldots + x_n = 0$), and  hence on $S^{n-2}$. The group $G = Q \times \A_n$  has
then a faithful action on the join 
$M^3 * S^{n-2} \cong S^{n+2}$, by joining the actions of $Q$ on $M^3$ and of
$\A_n$ on $S^{n-2}$ (with $Q$ acting trivially on $S^{n-2}$ and $\A_n$
trivially  on $M^3$).

\medskip

So $G$ admits a faithful, topological action on $S^{n+2}$ (which is not locally
linear). We will show that $G$ does not admit a faithful, linear action on 
$S^{n+2}$.  We fix a linear action of $G$ on $S^{n+2}$ and suppose, by
contradiction, that the action is faithful.

\medskip

The linear action of $G$ on $S^{n+2} \subset \R^{n+3}$  defines a
$(n+3)$-dimensional real representation of $G$. The induced linear
representation of $\A_n$ on
$\R^{n+3}$ splits as a direct sum of irreducible representations. Suppose first
that $n \ge 7$; then the only irreducible representations of $\A_n$ in
dimensions smaller than $n+3$ are the trivial representation and the standard
representation in dimension $n-1$, so there is an orthogonal decomposition
$\R^{n+3} = \R^4 \oplus \R^{n-1}$ where $\A_n$ acts trivially on the first
summand and by the standard representation on the second one.  The group $Q$
preserves this decomposition and commutes elementwise with $\A_n$.
Complexifying, we have a $G$-invariant decomposition 
$\C^{n+3} = \C^4 \oplus \C^{n-1}$; then, by Schur's Lemma, $Q$ acts by
homotheties on the second summand $\C^{n-1}$, i.e. by scalar multiples of the
identity (see [S] or [FH] for the representation theory of finite groups), and
hence the action of $Q$ on $\C^{n-1}$ factors through the action of an abelian
group. The abelianization of $Q$ is the Klein 4-group
$\Z_2 \times \Z_2$, generated by the images of $i$ and $j$; in particular, the
central involution $-1$ and the cyclic subgroup $\Z_{abc}$ of $Q$ act trivially
on $\C^{n-1}$, and hence also on  $\R^{n-1}$.

\medskip

Now we consider the action of $Q$ on the first summand $\R^4$. Since the action
of $Q$ on  $\R^{n+3} = \R^4 \oplus \R^{n-1}$ is faithful,  the central
involution $-1$ and the cyclic subgroup $\Z_{abc}$ of $Q$ act faithfully on the
first summand $\R^4$; then also
$i,j, k$ and  the subgroup  $Q(8)$ of $Q$ (the quaternion group of order 8) act
faithfully on $\R^4$, and hence the whole group $Q$. Since 
$Q$ is not isomorphic to a subgroup of ${\rm O}(4)$ (see [Mn]  and [Z3, section
3]), this is a contradiction, so the  action of $G$ on $S^{n+2}$ cannot be
faithful.

\medskip

We have proved the Theorem for $n \ge 7$, or $d \ge 9$. The cases $d = 8, 7$
and 6 are similar, considering the alternating and symmetric groups $\A_6$, 
$\S_5$ and $\A_5$ and their irreducible representations (see [Co] for their
character tables).  

\medskip

This concludes the proof of the Theorem.

\bigskip

{\bf Remarks.} i)  We discuss the case $d=5$. The Milnor group $Q$ considered
before admits a faithful action on $M^3 * S^1 \cong S^5$ (with the trivial
action on $S^1$). Now the authors of [CKS] remark (in a note in section 2) that
a Milnor group $Q(8a,b,c)$ is not isomorphic to a subgroup  of ${\rm O}(m)$,
for  $m \le 7$ (indicating an idea of a proof), so this would imply that $Q$
does not admit a faithful, linear action  on a sphere  of dimension less than
seven. The same holds then for the infinitely many groups  $Q \times \Z_k$; on
the other hand, these groups admit a faithful, topological action on $M^3 * S^1
\cong S^5$, by letting $\Z_k$ act faithfully by rotations on the second factor.

\medskip

ii) Elaborating on this, there is a faithful, topological  action of a group $G
= Q \times \Z_k \times \A_n$ on $M^3 * S^1 * S^{n-2} \cong S^5 * S^{n-2} \cong
S^{n+4}$. Now, if $Q$ is not a subgroup of ${\rm O}(6)$, exactly as in the
proof of the Theorem one shows that $G$ does not admit a faithful, linear
action on $S^{n+4}$, so in each dimension $d \ge 5$ there are infinitely many
groups as in the Theorem (considering the groups $Q \times \Z_k$ in dimensions 
$d<7$). 

\medskip

iii) In the proof of the theorem, the action of $Q$ on $M^3 * S^{n-2} \cong
S^{n+2}$, with fixed point set $S^{n+2}$, is not locally linear (otherwise
$Q$, which is not a subgroup of the orthogonal group ${\rm O}(4)$, would act
faithfully and orthogonally on a a 3-sphere orthogonal to the fixed point set). 
Choose a surjection $Q \to \Z_2$ and let a generator of $\Z_2$ act by minus
identity on $S^{n-2}$; this defines an action of $Q$ on $S^{n-2}$ without a
global fixed point. The kernel of the has fixed point set $S^{n-2}$ and is
isomorphic to a subgroup of ${\rm O}(4)$: is the combined action of $Q$ on 
$M^3 * S^{n-2}$ locally linear now (e.g., in the case
$M^3 * S^1 \cong S^5$)?  

\medskip

On the other hand, also the action of $\A_n$ on 
$M^3 * S^{n-2}$ is not locally linear: if an element of $\A_n$ has 0-dimensional
fixed point set $S^0$ in $S^{n-2}$ then its fixed point set in 
$M^3 * S^{n-2}$ is the suspension or double cone $M^3 * S^0$ which is a homology
manifold but not a manifold in the two cone points. So, in order to construct a
locally linear action of a group $Q \times A$ on $M^3 * S^{n-2}$,  one should
avoid such low-dimensional fixed point sets for an action of some group $A$ on
$S^{n-2}$.

\medskip

iv) The group $Q$ admits a free action on
$M^3*M^3$ which is homeomorphis to  $S^7$ (by the solution of the
higher-dimensional Poincar\'e conjecture): is this action conjugate to a free 
linear action of $Q$  on
$S^7$ (noting that $Q$ occurs as a fixed-point free subgroup of
${\rm SU}(4)$ and hence ${\rm O}(8)$)? Note that $Q$ has a free action also on
$S^{11} \cong S^7 * M^3 \cong M^3 * M^3 * M^3$ but no free linear action on
$S^{11}$.

\medskip

Finally, let $M^3$ be any homology 3-sphere with a free action of the cyclic
group $\Z_p$. Letting $\Z_p$ act by rotations on $S^1$, it has a free action on 
$M^3 * S^1 \cong S^5$: is this action conjugate to a linear action?

\bigskip 
\vfill \eject

\centerline {\bf References}

\medskip

\item {[BLP]}  M. Boileau, B. Leeb, J. Porti, {\it  Geometrization of
3-dimensional orbifolds.}  Ann. Math. 162,  195-250  (2005)

\item {[Ca]}  J.W. Cannon, {\it  The recognition problem: what is a topological
manifold.} Bull. Amer. Math. Soc. 84, 832-866  (1978)

\item {[CKS]}  W. Chen, S. Kwasik, R. Schultz, {\it Finite symmetries of
$S^4$.}  Forum Math. 28,  295-310  (2016)

\item {[Co]} J.H. Conway, R.T. Curtis, S.P. Norton, R.A. Parker, R.A.Wilson,
{\it Atlas of Finite Groups.} Oxford University Press 1985

\item {[D]}  M.W. Davis, {\it  A survey of results in higher dimensions.}  The
Smith Conjecture, edited by J.W. Morgan, H. Bass, Academic Press 1984,  227-240

\item {[DL]}  J. Dinkelbach, B. Leeb , {\it  Equivariant Ricci flow with
symmetry and applications to to finite group actions on 3-manifolds.}  Geom.
Top. 13, 1129-1173  (2009)

\item {[DoH]}  R.M. Dotzel, G.C. Hamrick, {\it  $p$-group actions on homology
spheres.} Invent. math. 62, 437-442  (1981)

\item {[FH]} W. Fulton, J. Harris, {\it Representation Theory: A First Course.} 
Graduate Texts in Mathematics 129,  Springer-Verlag 1991

\item {[GZ]} A. Guazzi, B. Zimmermann, \hskip 1mm {\it On finite simple groups
acting on homology spheres.}  Monatsh. Math. 169,  371-381 (2013)

\item {[MeZ1]} M. Mecchia, B. Zimmermann, {\it On finite simple and nonsolvable
groups acting on homology 4-spheres.} Top. Appl. 153,  2933-2942  (2006)

\item {[MeZ2]} M. Mecchia, B. Zimmermann, {\it On finite  groups acting on
homology 4-spheres and finite subgroups of ${\rm SO}(5)$.}  Top. Appl. 158, 
741-747  (2011)

\item {[Mg]} R.J. Milgram, {\it Evaluating the Swan finiteness obstruction for
finite groups.} Algebraic and Geometric Topology. Lecture Notes in Math. 1126,  
Springer 1985, 127-158

\item {[Mn]} J. Milnor, {\it Groups which act on $S^n$ without fixed points.}
Amer. J. Math. 79, 623-630  (1957)

\item {[Mu]} J.R. Munkres, {\it Elements of Algebraic Topology.} Addison-Wesley
Publishing Company 1984

\item {[S]}  J.-P. Serre, {\it Linear Representations of Finite Groups.}
Graduate Texts in Mathematics 42,  Springer 1977

\item {[Z1]} B. Zimmermann, {\it Some results and conjectures on finite groups
acting on homology spheres.} Sib. Electron. Math. Rep. 2, 233-238 (2005)  
(http://semr.math.nsc.ru)

\item {[Z2]} B. Zimmermann, {\it On the classification of finite groups acting
on homology 3-spheres.} Pacific J. Math. 217, 387-395 (2004)

\item {[Z3]} B. Zimmermann, {\it On finite simple groups acting on homology
spheres with small fixed point sets.}  Bol. Soc. Mat. Mex. 20,  611-621 (2014)

\bye